\documentclass[11pt,a4paper]{article}
\usepackage[utf8]{inputenc}

\title{Entropy Minimization for Optimization of Expensive, Unimodal Functions}
\author{Xiaohe Luo, Warren Powell}
\date{February 2023}
\usepackage{natbib}
\usepackage{setspace}
\usepackage{diagbox}
\usepackage{graphicx}
\usepackage{dsfont}
\usepackage{amssymb}
\usepackage{amsmath}
\usepackage{multirow}
\usepackage{multicol}
\numberwithin{equation}{section}
\usepackage{amsthm}
\usepackage{hyperref}
\usepackage{cleveref}
\usepackage{mathtools}
\usepackage{dsfont}
\usepackage{caption}
\usepackage{subcaption}
\usepackage{algorithmic}
\usepackage{algorithm}
\usepackage{geometry}
\usepackage[table,xcdraw, dvipsnmaes]{xcolor}
\newcommand{\E}{\ensuremath{\mathbb{E}}}

\newcommand{\X}{\ensuremath{\mathcal{X}}}
\newcommand{\Prob}{\ensuremath{\mathbb{P}}}

\newcommand\norm[1]{\left\lVert#1\right\rVert}

\newtheorem{theorem}{Theorem}[section]
\newtheorem{corollary}{Corollary}[theorem]

\newtheorem{lemma}[theorem]{Lemma}
\newtheorem*{remark}{Remark}
\theoremstyle{definition}
\newtheorem{definition}{Definition}[section]
\newtheorem{assumption}{Assumption}

\makeatletter
\newtheorem*{rep@theorem}{\rep@title}
\newcommand{\newreptheorem}[2]{%
\newenvironment{rep#1}[1]{%
 \def\rep@title{#2 \ref{##1}}%
 \begin{rep@theorem}}%
 {\end{rep@theorem}}}
\makeatother
\newreptheorem{theorem}{Theorem}

\begin{document}

\maketitle
\begin{abstract}
Maximization of an expensive, unimodal function under random observations has been an important problem in hyperparameter tuning. It features expensive function evaluations (which means small budgets) and a high level of noise.  We develop an algorithm based on entropy reduction of a probabilistic belief about the optimum. The algorithm provides an efficient way of estimating the computationally intractable surrogate objective in the general Entropy Search algorithm by leveraging a sampled belief model and designing a metric that measures the information value of any search point.
\end{abstract}

\section{Introduction}\label{intro} 
Optimization of noisy, unimodal functions which are expensive to compute has long been a hard problem to solve in the field of global optimization. The underlying problem can be formulated as follows:
\begin{equation}
\begin{aligned} \label{noiseless_obj}
    &\underset{x \in \mathcal{X} \subseteq \mathbb{R}}{\text{max}} && f(x) = \E[F(x,W)] \\
\end{aligned}
\end{equation}
where $\mathcal{X}$ is a bounded feasible domain and $W$ is a random variable representing the noise involved when observing the true function values. The black-box objective function $f$ is usually continuous, yet nonconvex, and its gradients are inaccessible. In this paper, we restrict our attention to finding the unique global optimum of this class of functions with a unimodal structure. In real-world applications, however, the measurement of the objective function can be both noisy and expensive: tuning hyperparameters of a machine-learning algorithm \citep{snoek2012practical}, arranging drug trials, controlling robots, and finding stepsizes in stochastic gradient algorithms, to name just a few applications. When faced with expensive functions, we have to solve these problems with a very limited number of observations. 

Bayesian optimization (BO) has been a popular strategy for solving the proposed problem due to its ability to find the near optimal solution in fewer experiments \citep{wu2016parallel,snoek2012practical}. BO achieves this goal by assuming a statistical model on the unknown, complicated objective function, which encodes the prior belief about the underlying function. This statistical model is updated after new observation(s) are made and reflects the incoming knowledge acquired as the experiment proceeds. In addition, Bayesian optimization adopts a comparatively tractable acquisition function as the new objective to guide the sampling process. 

At each iteration $n$, a BO algorithm seeks the point that maximizes the acquisition function as the next point for function evaluation, which means there is no attempt to maximize the value of information as is done with policies such as the knowledge gradient. The acquisition function can be identified as the objective function of the imbedded surrogate maximization problem used in three of the four classes of policies: cost function approximations (CFAs), value function approximations (VFAs) and direct lookahead policies (DLAs) \citep{powell2007approximate, powell2022unified}. 

Since \citet{mockus1994application} set up the theoretical foundation, Gaussian Process regression (GP) has become the typical statistical prior most BO algorithms adopt. The focus of the previous literature has been on designing an effective policy, equivalently the acquisition function, based on a GP prior \citep{williams2006gaussian}. Improvement-based policies proposed in previous works include probability of improvement (PI) \citep{kushner1964new}, expected improvement (GP-EI) \citep{mockus1978application, jones1998efficient} and knowledge gradient (KG) \citep{scott2011correlated, frazier2009knowledge}. These are direct lookahead policies (DLAs in \cite{powell2022unified}) as they use the one-step improvement in the estimation of the global optimum, measured by some metric, as the criterion to pick the next query point. Policies based on other criteria are, for example, upper confidence bound (GP-UCB) \citep{srinivas2009gaussian}, which is a popular method in the CFA class and entropy search (ES) \citep{hennig2012entropy}. Entropy search is particularly interesting because it hybridizes the methodology of one-step lookahead and the idea of finding a good substitute for the unknown value functions. Each of these policies incorporates a distinct philosophy of sampling a good sequence of queries that help one decide on the final optimum. 
While all of these methods aim at finding the next point for observation, the resulting surrogate objective functions are, though better than the original objective function $f$, not always easy to maximize. For example, KG and ES involve maximizing the expectation of a complicated, nonconvex function. Let $\mathcal{D}^n \coloneqq \{x^i,\hat{f}(x^i)\}_{i=0}^{n-1}$ be the set of data available at iteration $n$. Specifically, the surrogate objective function adopted by ES takes the form:
\begin{equation}
\begin{aligned} \label{ES}
    ES^n(x) =  H(p(X^*|\mathcal{D}^n)) - \E_{\hat{f}(x)|\mathcal{D}^n}[H(p(X^*|\mathcal{D}^n \cup \{x,\hat{f}(x)\}], \\
\end{aligned}
\end{equation}
where $\hat{f}(x)$ is a noisy observation at $x$, $p(X^*|\mathcal{D}^n)$ is the posterior distribution of the global optimizer $x^*$ at iteration $n$ and $H(\cdot)$ is the Shannon differential entropy \citep{hennig2012entropy}. Formula \eqref{ES} does not have an analytical expression if the GP prior is selected, resulting in a computationally intractable objective \citep{frazier2018tutorial,hernandez2014predictive}. As a consequence, it becomes numerically expensive to estimate \eqref{ES} not only because Monte Carlo sampling is required in optimizing this kind of surrogate objective functions \citep{brochu2010tutorial} but also because each observation is hard to compute. To address the computational difficulty of ES, there have been several attempts in the literature. By recognizing \eqref{ES} to be mutual information, \citet{hernandez2014predictive} proposes Predictive Entropy Search (PES), whose surrogate objective function is theoretically identical to \eqref{ES} by the property of mutual information: 
\begin{equation}
\begin{aligned} \label{PES}
    PES^n(x) =  H(p(\hat{f}(x)|\mathcal{D}^n)) - \E_{X^*|\mathcal{D}^n}[H(p(\hat{f}(x)|\mathcal{D}^n, X^*)]. \\
\end{aligned}
\end{equation}
A further step is taken by considering the maximum function value $y^*$ instead of $x^*$ in formula \eqref{PES} \citep{hoffman2015output, wang2017max}, which yields the Max-value Entropy Search (MES):
\begin{equation}
\begin{aligned} \label{MES}
    MES^n(x) =  H(p(\hat{f}(x)|\mathcal{D}^n)) - \E_{Y^*|\mathcal{D}^n}[H(p(\hat{f}(x)|\mathcal{D}^n, Y^*)]. \\
\end{aligned}
\end{equation}
Even though both of the above methods simplify the procedure of approximating \eqref{ES}, equations \eqref{PES} and \eqref{MES} still require Monte Carlo methods to estimate the expectation in the formulas as well as additional techniques such as expectation propagation to sample from the distribution $p(X^*|\mathcal{D}^n)$ or  $p(Y^*|\mathcal{D}^n)$.

In this paper, we present a new BO algorithm based on the idea of entropy search: the sampled-belief entropy search (SBES). By assuming a sampled belief model on the underlying function $f$, we convert the hard-to-evaluate ES surrogate objective function into a deterministic function that can be optimized via well-established deterministic, derivative-free optimization methods. This conversion is established on the elimination logic of Fibonacci search \citep{ferguson1960} and the probability of observing the correct gradient signs introduced for a gradient-based method in \citet{powell2012optimal, peter2013bisection}. In the derivative-free setting, we extend the probability of observing correct gradient signs to the probability of observing the correct location of the optimum, also called the probability of correct region assignment. We show that by carefully designing the probability of correct region assignment, SBES outperforms when the total budget is small.

This paper makes the following contributions: 1)We design a new entropy-search based policy that specifically tackles the difficulty in finding the optimum of unimodal, noisy and expensive functions. To combat the challenge introduced by noisy observations and to fully leverage the proposed unimodal structure, a discrete parametric prior is adopted to model the truth function $f$. 2)By introducing the probability of correct region assignment that can be calculated using the parametric prior, we derive an updating rule for any posterior $p(X^*|\mathcal{D}^n)$ and thus successfully turn expression \eqref{ES} into an analytical formula. 3)We present an error bound on the one-step information gain of SBES in the stochastic setting. 4) We conduct empirical experiments to compare SBES and other BO algorithms, including both the non entropy-search based and the entropy-search based methods, under different truth functions across various levels of noise. We show that SBES is robust,  competitive with other BO algorithms at high noise levels and outperforms at low and medium noise levels.

The paper is organized as follows. Section 2 describes the process of solving problem \eqref{noiseless_obj} formally as a one-dimensional search problem and introduces our proposed models on this problem. The modeling part includes the definition of the probability of correct region assignment function, the sampled belief model and the probability distribution over the belief of $x^*$. Section 3 starts with formulating the one-dimensional search problem as a stochastic sequential decision problem. It is then followed by a discussion of the surrogate objective function we use and ways of solving it when (i) $\mathcal{X}$ is finite and discrete, or (ii) $\mathcal{X}$ is compact. Section 4 is devoted to establishing the one-step error bound of our algorithm. This paper concludes with section \ref{section5:experiment} in which we compare SBES against popular bench-mark algorithms such as GP-UCB, GP-EI and the response surface method.

\section{Problem and Models}\label{prob_model} 
We begin by defining the problem of finding the location of the optimum of a unimodal function mathematically. We formulate the problem of finding the best algorithm (policy) $\pi$ for learning the value $x^{\pi,N}$ that maximizes $\E [F(x,W)]$ in a budget of $N$ function evaluations. The second part of this section is devoted to the statistical model we use to approximate the underlying function $f$ and the modeling of uncertainty based on this statistical model. 

\subsection{Problem Definition.}\label{problem_def} 
Suppose there is a unimodal function $f:\mathbb{R}\rightarrow\mathbb{R}$ with a feasible domain $\mathcal{X}$ that can be finite and discrete or compact. There are only noisy observations of this function available whenever we decide to evaluate the function at some points. Let $\hat{f}(x)$ be the observation of the function $f$ at a point $x \in \mathcal{X}$. Mathematically, define $F(x,W) \coloneqq \hat{f}(x)$, where $W $ is a random variable that inherits the randomness of the noise. We then make the following assumptions about $F(x,W)$:
\begin{enumerate}
    \item[(1)]Unbiasedness: at any point $x\in \mathcal{X}$, $\mathbb{E}[\hat{f}(x)]=f(x)$.
    \item[(2)]Homoscedasticity: for any $x$ and $y$ $\in \mathcal{X}$, $Var(\hat{f}(x))=Var(\hat{f}(y))=\sigma^2$.
    \item[(3)]Gaussian noise: $\hat{f}(x) = f(x) + \epsilon$, where $\epsilon \sim \mathcal{N}(0,\sigma^{2})$.
    \item[(4)]$f(x)$ is continuous and might be differentiable but the gradient of $f(x)$ is not available to us.
\end{enumerate}

Define $x^* \coloneqq \mathrm{argmax}_{x \in \mathcal{X}} \mathbb{E}[\hat{f}(x)]$. Our challenge is to find the optimum of this unimodal function $x^*$ with a limited number of measurements. Assume there is an algorithm $\pi$ that produces a solution $x^{\pi,N}$ after $N$ iterations. Then our goal is to find the best algorithm that solves:
\begin{equation} \label{final_reward}
\begin{aligned}
   \max_\pi \mathbb{E} \left\{F(x^{\pi,N},\hat{W}) |S^0\right\}.
\end{aligned}
\end{equation}
given an initial state $S^0$. In order to determine $x^{\pi, N}$ after $N$ iterations, the policy $\pi$ needs to do learning via exploration in the search region. Ideally, the optimal policy $\pi$ already knows the location of $x^*$ after $N$ experiments so that it will map $x^{\pi, N}$ to $x^*$. With this being said, $\pi$ also needs to determine the set of points for function evaluation, denoted as $\left \{ x^n \right \}_{n=1}^{N-1}$, up to iteration $N$ so that it can give the best estimate of $x^*$. Hence, taking the uncertainty in the initialization $S^0$ into account, problem \eqref{final_reward} is equivalent to:
\begin{equation} \label{expand_final_reward}
\begin{aligned}
   \max_\pi \mathbb{E}_{S^0}\mathbb{E}_{W^1,...,W^N|S^0}\left\{ \E_{\hat{W}} [F(x^{\pi,N},\hat{W}) |S^0 ]\right\}.
\end{aligned}
\end{equation}

\subsection{Models.}\label{model} 
In Fibonacci search, two function evaluations can provide an indication of where the optimum might be and thus eliminate a section of the region in $\X$ that the optimum is not in. This is a property that comes from unimodality and the assumption of no noise. Starting with two initial points, Fibonacci search chooses one point at each iteration in the search region to evaluate. Then this new function evaluation is used along with the previous function evaluation to narrow the region where $x^*$ might be located. 

In the stochastic setting, we still want to exploit the property of unimodality and extract the information about the location of $x^*$ from two function evaluations. However, the presence of noise hinders the elimination of regions as is done with classical Fibonacci search, so we maintain a dynamic belief about $x^*$ in the form of distribution instead. This belief is updated by choosing one point from the history of prior observations, and another point $z^n$ at which we perform another function evaluation.  The point $z^n$ is chosen to minimize the expected entropy in our belief about the location of $x^*$.

For the purpose of this section, we defer the discussion of how exactly we choose those points at each iteration $n$ to the section \ref{algo}. For now, assume that at iteration $n$, we pick a new point $z^n \in \mathcal{X}$ and a point from history $h^n \in H^n$, where $H^n \coloneqq   \{x^0 \} \cup \{z^i\}_{i=1}^{n-1}$ is the set of historically chosen points up to time $n$ and $x^0$ is the set of initial points. So the set of points we pick at iteration $n$ is: $x^n = (h^n, z^n)$. We then perform an (expensive) function evaluation $\hat{f}(z^n) = F(z^n, W^{n+1})$, and define the history of observations to be $\hat{f}^n \coloneqq  \{\hat{f}(x^0) \} \cup \{\hat{f}(z^i)\}_{i=1}^{n-1}$. A prior distribution of belief of the location of $x^*$ can be constructed using $(H^n, \hat{f}^n)$, denoted as $P^n$. The posterior distribution $P^{n+1}$ can then be calculated given the comparison between $\hat{f}(h^n)$ and $\hat{f}(z^n)$ using Bayes theorem.
In this section, the relative location of $h^n$ and $z^n$ is important; that is, whether $h^n < z^n$ or $h^n > z^n$ determines how we update $P^n$. So we label the smaller of the two points $x^n_l$ and the larger of the two points $x^n_r$. Then $x^n=(h^n, z^n)=(x^n_l, x^n_r)$.

In the remaining sections, we first introduce the sampled belief model used to represent the ground truth $f(x)$. Then, we discuss two probabilistic approaches that are used to model the uncertainty in function observations and the location of $x^*$. 

\subsubsection{Sampled Belief Model.}\label{sample_belief} 
\begin{figure}
    \centering
    \includegraphics[width=10cm]{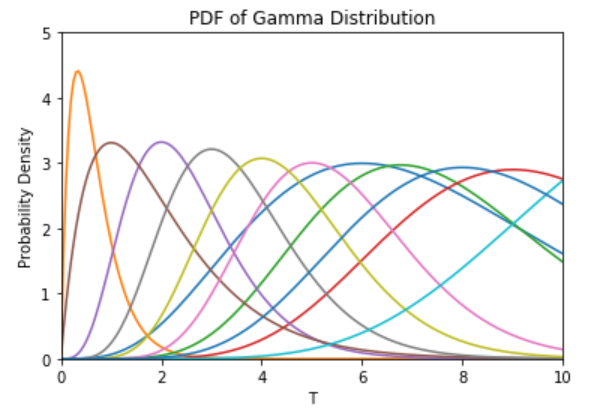}
    \caption{A family of gamma functions. \label{fig:gamma}}
\end{figure}

To address the question of how to estimate the true function appropriately, one approach is to use a finite set of parametric, unimodal functions (see figure \ref{fig:gamma} for a family of gamma distributions as an example). Suppose the ground truth can be parameterized as $f(x | \theta)$. Let $F_\Theta = \left \{ f_{k} = f(x |\theta_k) : \theta_k \in \Theta, \forall \ 1 \leq k \leq K\right \}$ be a parametric family of unimodal functions that are used to approximate $f(x |\theta)$. Also define $p^n_k$ to be the probability that $f_{k}$ is the best representation of $f$ at iteration $n$: 
\begin{align} \label{def:prob_theta}
    p^n_k \coloneqq \Prob[\theta = \theta_k|\mathcal{D}^n].
\end{align}
We then approximate $f$ by $\bar{f}^n(x) = \sum_{k=1}^K p_k^n f_k(x)$ for all $x \in \X$. 

During the experiment, $\bar{f}$ will become more precise as more data points come in. Initialize $p^{0}_k=\frac{1}{K}$, and implement Bayesian updating to $\{p^n_k \}$ after each observation $\hat{f}(z^n)$. That is:
\begin{align} \label{update_sample_belief}
     p_k^{n+1} = \frac{p(\hat{f}(x)|x=z^n,\theta_k)p_k^n}{\sum_{k=1}^K p(\hat{f}(x)|x=z^n,\theta_k)p_k^n} \ \ \forall 1 \leq k \leq K,
 \end{align}
 where $p(\hat{f}(x)|x=z^n,\theta_k)$ is the density function of a normally distributed random variable $\hat{f}(z^n) \sim \mathcal{N}(f_{k}(z^n), \sigma^2)$.

\subsubsection{Probability of Correct Region Assignment.}\label{prob_assign} 
Since the underlying function is unimodal, two function evaluations can provide us information about the location of $x^*$. Nevertheless, the noise in observations distorts this piece of information, so it is necessary to model the uncertainty in comparing two observations. Given any two points $x, y$ and their corresponding observations of $f$: $\hat{f}(x), \hat{f}(y)$, there are two possible outcomes: either $\hat{f}(x)>\hat{f}(y)$ or $\hat{f}(x)\leq \hat{f}(y)$. Similarly, the true function values evaluated at these two points can also be split into two cases: $f(x)>f(y)$ or $f(x)\leq f(y)$. Let $g(x,y)$ denote the probability that we are able to observe the true comparative relation between $f(x)$ and $f(y)$, i.e. if $f(x) > f(y)$, the probability of observing $\hat{f}(x) > \hat{f}(y)$ is $g(x,y)$; if $f(x)\leq f(y)$, then the probability of observing $\hat{f}(x) \leq \hat{f}(y)$ is $g(x,y)$.
Mathematically for any $x,y \in \mathcal{X}$,
\begin{align} \label{eq:g_fun}
g(x,y) \coloneqq \left\{\begin{matrix}
\mathbb{P}(\hat{f}(x)>\hat{f}(y)) &\quad \text{if} \ f(x)>f(y) \\ 
\mathbb{P}(\hat{f}(x) \leq \hat{f}(y)) &\quad \text{if} \ f(x) \leq f(y).
\end{matrix}\right.
\end{align}
The following truth table describes the relationships between $\hat{f}(x), \hat{f}(y)$ and $g(x,y)$:
\renewcommand{\arraystretch}{1.5}
\begin{center}
\begin{tabular}{ c|c|c } 
 \diagbox[width=10em]{Observation}{Truth}& $f(x)>f(y)$ & $f(x) \leq f(y)$ \\ 
 \hline
 $\hat{f}(x)>\hat{f}(y)$ & $g(x,y)$ & $1-g(x,y)$ \\ 
 \hline
 $\hat{f}(x)\leq \hat{f}(y)$ & $1-g(x,y)$ & $g(x,y)$ \\ 
\end{tabular}
\end{center}
Note that the function $g$ preserves symmetry by definition: $g(x,y) = g(y,x)$.

Recall that the problem of interest is where $x^*$ resides in the domain $\mathcal{X}$. Without loss of generality, let $x_l \coloneqq \text{min} \left \{x,y \right\}$ and $x_r\coloneqq \text{max} \left \{x,y \right\}$. By unimodality, $f(x_l)>f(x_r)$ implies $x^* < x_r$; $f(x_l) \leq f(x_r)$ implies $x_l \leq x^*$. However, $\left \{ x^* < x_r \right \}$ and $\left \{ x_l \leq x^* \right \}$ are not complementary events since there is a possibility that $x_l \leq x^* < x_r$. To define the conditional probabilities dependent on events of $x^*$, considering only the two cases,\ $f(x_l)>f(x_r)$ and $f(x_l)\leq f(x_r)$, is not sufficient. Instead, for any points $x,y$, define 
\begin{align} \label{eq:2}
    \bar{g}(x,y)
    \coloneqq \left\{\begin{matrix}
\mathbb{P}(\hat{f}(x_l)>\hat{f}(x_r)) &\quad \text{if} \ x_l <x^*< x_r\\ 
\frac{1}{2} &\quad \text{if} \ x^* \notin (x_l,x_r)
\end{matrix}\right.
\end{align}
 and consider the following table:
\renewcommand{\arraystretch}{1.5}
\begin{center}
\begin{tabular}{ c|c|c|c } 
 \diagbox[width=10em]{Observation}{Truth}& $x^* \leq x_l$ & $x_l <x^*< x_r$ & $x_r \leq x^*$\\ 
 \hline
 $\hat{f}(x_l)>\hat{f}(x_r)$ & $g(x_l,x_r)$ & $\bar{g}(x_l,x_r) $ & $1-g(x_l,x_r)$ \\ 
 \hline
 $\hat{f}(x_l)\leq \hat{f}(x_r)$ & $1-g(x_l,x_r)$ & $1-\bar{g}(x_l,x_r) $ & $g(x_l,x_r)$ \\ 
\end{tabular}
\end{center}
 The above truth table is constructed by assigning $\mathbb{P}(\hat{f}(x_l)>\hat{f}(x_r)|x^* \leq x_l) = g(x_l,x_r)$, $\mathbb{P}(\hat{f}(x_l) \leq \hat{f}(x_r)|x_r \leq x^*) = g(x_l,x_r)$ and $\mathbb{P}(\hat{f}(x_l) > \hat{f}(x_r)|x_r \leq x^* \leq x_r) = \bar{g}(x_l,x_r)$. Note that the equalities here don't hold mathematically but rather an assignment. We are using $g$ and $\bar{g}$ to approximate the true, desired probabilities. We can do so because the event $\left \{x^* \leq x_l \right \}$ implies $ \left \{f(x_l) \leq f(x_r) \right \}$ and $\left \{x_r \leq x^* \right \}$ implies $ \left \{f(x_r) \leq f(x_l) \right \}$.
 
\subsubsection{Probability Distribution over Belief of the Location of $x^*$.}\label{prob_x} 
The randomness in the location of the optimum is induced by not only the noise in observations but also the randomness in the function space to which $f$ belongs to. Suppose there is a set of unimodal functions $\mathcal{F}$ from which $f$ is generated. Then $\mathcal{F}$ is naturally endowed with a measure based on prior knowledge of the likelihood of all candidate functions. Moreover, this distribution can be updated after receiving data $\mathcal{D}^n$. Following the same line of logic, let $X^*$ be a random variable corresponding to our belief about the real optimum $x^*$. The probability measure $P^n$, which encodes our knowledge of $x^*$ given historical data $(H^n,\hat{f}^n)$ up to iteration $n$, is thus well-defined as introduced in \citep{hennig2012entropy}. Note that given a fixed $x^*$, the true density of $X^*$ over the domain $\mathcal{X}$ is the Dirac delta function:
\begin{align} \label{true_pdf}
p_{X^*}(x | x^*) = \delta(x-x^*) 
\end{align}
(if $\mathcal{X}$ is finite, $p_{X^*}(x|x^*) = \mathds{1}_{\{x=x^*\}}$). 
Let $S^n = (H^n,h^n,P^n)$ be the random variable including information from the history up through $n$. Then define $P^{n+1}(x)$ to be the probability distribution of $X^*$ based on the history of observations and comparisons, the decision we make $x^n$ and the randomness in the $n+1^{th}$ experiment $W^{n+1}$. Mathematically,
\begin{align} \label{def_posterior}
    P^{n+1}(x) dx = \mathbb{P}(X^* \in dx | S^n, x^n, W^{n+1}).
\end{align}
 At each iteration $n$, we are trying to approximate $p_{X^*}$ by $P^n$ given observations. Equation \eqref{def_posterior} indicates that at iteration $n$, after we have made the decision $x^n$ and have the two function evaluations we need $\hat{f}(x^n)=(\hat{f}(h^n), \hat{f}(z^n))$, $P^n$ can be updated to $P^{n+1}$ by Bayes' theorem. In the following, we  present the updating rule in detail.
 
 First we define an auxiliary variable $\hat{y}^{n+1}$ as:
 \begin{align} \label{def_yhat}
     \hat{y}^{n+1} = \mathds{1}_{\left \{ \hat{f}(x_l^n) \leq \hat{f}(x_r^n)  \right \}}.
 \end{align}
 Note that $\hat{y}^{n+1}$ is a function of both the decision $x^n$ and the randomness $W^{n+1}$, indicating the comparative relation between the two observations that are used to update our belief. Given the definition of $\hat{y}^{n+1}$, we can rewrite \eqref{def_posterior} as:
 \begin{align} \label{formula_posterior}
    P^{n+1}(x)dx 
    & = \mathbb{P}(X^*\in dx |S^n, x^n, \hat{y}^{n+1}) \nonumber \\
    &=\frac{\mathbb{P}(\hat{y}^{n+1}, X^* \in dx | S^n, x^n)}{\mathbb{P}(\hat{y}^{n+1} |  S^n, x^n)} \nonumber \\
    &=\frac{\mathbb{P}(\hat{y}^{n+1} | X^*\in dx,  S^n, x^n) P^n(x)}{\int_{\mathcal{X}} \mathbb{P}(\hat{y}^{n+1} | X^* \in dx,  S^n, x^n) P^n(x)dx}  dx .
\end{align}
Equation \eqref{formula_posterior} is a general formula for calculating the posterior $P^{n+1}$ given a prior $P^n$. With the help of the sampled belief model and the probability of correct region assignment, \eqref{formula_posterior} has an analytic expression (please see appendix for detailed calculations). The exact updating formula can be found in section \ref{five_elements}.

\section{Algorithm and Policy} \label{algo} 
With the models introduced in Section 2, it is possible to solve problem \eqref{expand_final_reward} with a well-defined algorithm. We first formalize the SBES algorithm in the beginning of this section. Then, we present all the necessary components, such as the updating rules for all state variables, of the algorithm with the following steps: begin with defining the five basic elements of the sequential decision problem \eqref{expand_final_reward}; then design a one-step lookahead policy $\pi$ for solving this sequential decision problem. We then show that the policy itself is generated by a surrogate optimization problem.
\subsection{Algorithm.}
Section 2 gives a brief overview of the SBES algorithm and all the probabilistic models we use to quantify the uncertainty in our problem without specifying how we make the decision $x^n$. We now describe the full algorithm along with the surrogate objective function. Inspired by Fibonacci search, our algorithm starts with a pair of initial points $x^0=(x^0_l, x^0_r)$ and pick a new point in the search region, denoted as $z^n$ at iteration $n \geq 1$ to obtain a new function evaluation $\hat{f}(z^n)$. This new observation is compared with one of the previous observations to update the belief $P^n$. This means that besides $z^n$, we also pick a historical point $h^n \in H^{n}$ and use the corresponding function evaluation $\hat{f}(h^n) \in \hat{f}^{n}$ directly. These two points are chosen to minimize the expected entropy of our belief $P^{n+1}$ in the next step. See algorithm \ref{SBES} for details.
\begin{algorithm}
   \caption{Sampled Belief Entropy Search}
   \label{SBES}
\textbf{Initializations:} draw an initial set of samples $x_0 = (x_l^0,x_r^0)$ in $\mathcal{X}$; pick an initial state $S^0 = \left (P^0, H^0, \hat{f}^{0}, g^0, \bar{g}^0, \left \{p_k^0 \right \}_{k=1}^K \right)$.
\begin{algorithmic}[1]
   \FOR{$n=1$ {\bfseries to} $N$}
    \IF{n=1} \STATE {Evaluate $x^0$ and obtain $\hat{f}(x^0)=(\hat{f}(x_l^0), \hat{f}(x_r^0))$}.
    \ELSE \STATE{Evaluate $z^{n-1}$ and observe $\hat{f}(z^{n-1})$}; obtain $\hat{f}(h^{n-1})$ from $\hat{f}^{n}$. \ENDIF
   \STATE Update $S^{n-1}$ to $S^n$ according to equations (18)-(26).
   \STATE Sample $m$ points from $P^n$ to form the set $\mathcal{A}^n$. Pick $x^{n} = (h^{n}, z^{n}) = \mathop{argmin}_{h \in H^n,  z\in \mathcal{A}^n} \nu^n(h,z)$ defined in \eqref{eq:sbes}.
   \ENDFOR
\end{algorithmic}
\textbf{Return} $\bar{x}^N = \mathop{argmax}_{x \in \mathcal{X}} P^N(x)$.
\end{algorithm}

\subsection{The Sequential Learning Process.} \label{five_elements} 
Our algorithm is a form of sequential decision problem, which can be described by five core elements (see \citep{powell2022unified}): state variables, decision variables, exogenous information, transition function and objective function.  We describe these below:
\begin{description}
\item[\textbf{State Variables}] $S^n = \left (P^n, H^n, \hat{f}^{n}, g^n, \bar{g}^n, \left \{p_k^n \right \}_{k=1}^K \right)$.\\
At iteration $n$, the state variable $S^n$ consists of two parts. The first part is our belief about some important quantities in our sequential decision problem. This includes: the belief $P^n$ of the location of the optimum $x^*$ given the first n observations, the belief about the probabilities of correct region assignment $(g^n, \bar{g}^n)$, and the belief about how ``close" each sample curve is to the underlying function $\left \{p_k^n \right \}_{k=1}^K$. The second part of the state variables involves the history of our experiments: $H^n$ denotes the set of points we have chosen up to iteration $n$ and $\hat{f}^n$ is the corresponding observations of the function value at those points. This is all the information we need to make a decision at iteration $n$.

\item[\textbf{Decision Variables}] $x^n = (h^n, z^n)$.\\ 
Due to the logic that we only choose one new point for function evaluation and the other function evaluation comes from the history, the decision variable has two components. It is straightforward that the first decision variable is $z^n \in \X$, the point we pick for function evaluation. The second decision is the observation $h^n \in H^n$ we use from history. After $n$ observations, we make the decision of the next point to observe using our policy $x^n = X^\pi(S^n)$ which depends on the information in $S^n$. 

\item[\textbf{Exogenous Information}] $W^{n+1}$.\\
The variable $W^{n+1}$ is the new information that is observed from the ${n+1}^\text{st}$ function evaluation which means we can write our sequential decision process as: $(S^0, x^0, W^1, S^1, x^1, ...,$ $ S^{N-1}, x^{N-1},W^N, S^N, x^N)$. After being in a state $S^{n}$ and choosing an action $x^{n}$, we will observe a realization of the random variable $W^{n+1}$ coming outside of the system. Since the noise in a single observation follows $\epsilon \sim \mathcal{N}(0,\sigma^{2})$ and when $n=0$ we have two new observations, $W^{1}$ is a two-dimensional Gaussian variable such that $W^{1} \sim \mathcal{N}(0, \sigma^2 I)$. When $n \geq 1$, $W^{n+1}$ is the randomness in evaluating $\hat{f}(z^n)$. So $W^{n+1} \sim \mathcal{N}(0, \sigma^2).$ 

\item[\textbf{Transition Function}] $S^{n+1} = S^M(S^n, x^n, W^{n+1})$.\\
The transition function describes how a state $S^n$ evolves to $S^{n+1}$ given the decision $x^n$ and the exogenous information $W^{n+1}$. At iteration $n$, the state variables in this problem is $S^n =$ ($P^n, H^n, h^n, g^n, \bar{g}^n, \left \{p_k^n \right \}_{k=1}^K$). This means that the transition function consists of a series of equations for updating each state variable, as follows:
\begin{enumerate}
    \item  \textbf{$P^n$ to $P^{n+1}$}: following the arguments in section \ref{prob_x} and letting $F^n$ be the CDF of $P^n$, \eqref{formula_posterior} takes the forms:\\
    \textbf{case 1:} if $\hat{y}^{n+1} = 1$,
    \begin{align} \label{eq:3.1}
        P^{n+1}(x) &=
        \left\{\begin{matrix}
        \frac{1-g^n(x_l^n,x_r^n)}{U_1^n(x_l^n,x_r^n)}P^n(x) 
        &\text{if} \ x \leq x_l^n,\\
        \frac{1-\bar{g}^n(x_l^n,x_r^n)}{U_1^n(x_l^n,x_r^n)}P^n(x)  & \quad \text{if} \ x_l^n<x<x_r^n, \\
        \frac{g^n(x_l^n,x_r^n)}{U_1^n(x_l^n,x_r^n)}P^n(x) 
        &\text{if} \ x_r^n \leq x.
        \end{matrix}\right.
    \end{align}
    \begin{align} \label{eq:u1}
U_1^n(x_l^n,x_r^n) \coloneqq 
& (1-g(x_l^n,x_r^n))F^n(x_l^n) +(1-\bar{g}(x_l^n,x_r^n))(F^n(x_r^n)-F^n(x_l^n)) \nonumber \\
 & + g(x_l^n,x_r^n)(1-F^n(x_r^n)).
\end{align}
    \textbf{case 2:} if $\hat{y}^{n+1} = 0$,
    \begin{align} \label{eq:3.2}
        P^{n+1}(x) &=
        \left\{\begin{matrix}
        \frac{g^n(x_l^n,x_r^n)}{U_0^n(x_l^n,x_r^n)}P^n(x) & \text{if} \ x \leq x_l^n,\\ 
        \frac{\bar{g}^n(x_l^n,x_r^n)}{U_0^n(x_l^n,x_r^n)}\ P^n(x) &\quad \text{if} \ x_l^n<x<x_r^n, \\
        \frac{1-g^n(x_l^n,x_r^n)}{U_0^n(x_l^n,x_r^n)} \ P^n(x) & \text{if} \ x_r^n \leq x. \end{matrix}\right. 
    \end{align}
    \begin{align}\label{eq:u0}
U_0^n(x_l^n,x_r^n) \coloneqq
& g(x_l^n,x_r^n)F^n(x_l^n) + \bar{g}(x_l^n,x_r^n)(F^n(x_r^n)-F^n(x_l^n)) \nonumber \\
& + (1-g(x_l^n,x_r^n))(1-F^n(x_r^n)).
\end{align}
    \item \textbf{$H^n, h^n$ to $H^{n+1}, h^{n+1}$}:
    \begin{align}
        H^{n+1} = H^n \cup z^n. \nonumber\\
        h^{n+1} = h^n \cup \hat{f}(z^n). \nonumber
    \end{align}
    \item \textbf{$g^n, \bar{g}^n$ to $g^{n+1}, \bar{g}^{n+1}$}: $\forall x,y \in \X$,
    \begin{align}
        & g^{n+1}(x,y) = \sum_k^K p^{n}_k \cdot g_k(x,y). \nonumber\\
        & \bar{g}^{n+1}(x,y) = \sum_{k=1}^K  \mathbb{P}(\theta_k | x<x^*<y) \left(1-\Phi(-\frac{f_k(x)-f_k(y)}{\sqrt{2}\sigma})\right). \nonumber
    \end{align}
    where,
    \begin{align}
        g_k(x,y)=
        \mathds{1}_{\left \{f_k(x) \geq f_k(y)\right \}} \left(1-\Phi(-\frac{f_k(x)-f_k(y)}{\sqrt{2}\sigma})\right) + \mathds{1}_{\left \{f_k(x) < f_k(y)\right \}} \Phi(-\frac{f_k(x)-f_k(y)}{\sqrt{2}\sigma}). \nonumber
    \end{align}
    \begin{align}
       & \bar{\Theta} \coloneqq \left \{\theta_k: x_k^* \in (x, y),\  \forall 1 \leq k \leq K \right \}. \nonumber \\
            &\mathbb{P}(\theta_k | x<x^*<y) = \left\{\begin{matrix}
            0 & \text{if} \ \theta_k \notin \bar{\Theta},\\ 
            \frac{p_k^{n}}{\sum_{\theta_k \in \bar{\Theta}} p_k^{n}} & \text{if} \ \theta_k \in \bar{\Theta}.
        \end{matrix}\right. \nonumber
    \end{align}
    \item \textbf{$\left \{p_k^n \right \}_{k=1}^K$ to $\left \{p_k^{n+1} \right \}_{k=1}^K$}: please refer to \eqref{update_sample_belief}.
\end{enumerate}

\item[\textbf{Objective Function.}] The objective of this problem is to find a policy $\pi^*$ that solves:
\begin{equation}
\begin{aligned}
   \max_\pi \mathbb{E} \left\{F(x^{\pi,N},\hat{W}) |S^0\right\}. \nonumber
\end{aligned}
\end{equation}
\end{description}

\subsection{Policy in the SBES Algorithm.} \label{sbes_policy} 
We now present the policy for choosing the point $x^n$. We first describe the logic behind Entropy Search (ES) and the surrogate optimization problem that optimizes the expected entropy reduction, followed by a discussion of how to solve this surrogate optimization problem when $\mathcal{X}$ is a finite set or a compact interval.

\subsubsection{One-step Lookahead Entropy Reduction.} 
Recall that in the deterministic case, the metric we use for determining the distance between the true optimum $x^*$ and our estimate of the optimum $x^N$ is the $\mathbb{L}^1$ norm: $|x^N-x^*|$. In the stochastic setting, we want a similar ``distance" metric on how close our estimated density function $P^n$ is to the true density of $X^*$: $p_{X^*}$. However, without knowing the exact location of $x^*$, there is no way of comparing $P^n$ and $p_{X^*}$. To address this problem, \citet{hennig2012entropy} suggest using the Kullback-Leibler(KL) divergence, also called \emph{relative entropy}, from the uniform measure $U_{\mathcal{X}}$ to the posterior $P^n$ as a mean for assessing the current information about $x^*$.

Denote the differential entropy of any density function $p$ on $X^*$ by 
\begin{equation}
\begin{aligned}
H(p)= -\int_{\mathcal{X}}p(x)log_2(p(x))dx.
\end{aligned}
\end{equation}
The relative entropy (KL-divergence) from $q$ to $p$ is defined as:
\begin{equation}
\begin{aligned}
D_{KL}(p\parallel q) &= \int_{\mathcal{X}}p(x)log_2(\frac{p(x)}{q(x)})dx \nonumber \\
                     &= -\int_{\mathcal{X}}p(x)log_2(q(x))dx - H(p).
\end{aligned}
\end{equation}
Relative entropy is a measure of how ``close" two distributions are. $D_{KL}(p \parallel q) = 0$ if and only if $p(x) = q(x)$ almost everywhere. On the other hand, the relative entropy from the uniform distribution $U_{\mathcal{X}}$ to the dirac delata function $\delta(x-x^*)$ is defined to be $\infty$. Hence, a larger relative entropy from the uniform measure to $P^n$ implies greater dissimilarity between them, meaning we have more information about the location of $x^*$. This inspires us to use the expected relative-entropy maximization as the value function of the surrogate problem.

With relative entropy as the metric, we can now define the contribution functions $\Tilde{C}_i$ for our surrogate problem. Recall that for any $1 \leq n \leq {N}$, the value function of a state $s$ is $\Tilde{V}_n^\pi(s) = \mathbb{E}[\sum_{i=n}^{N} \Tilde{C}_i(S^i, X^\pi(S^i), W^{i+1}) | S^n=s, x^n=X^\pi(S^n)]$. For any $n \leq i \leq N$ and the state variable $S^i$, define 
\begin{align}
    \tilde{C_i}(S^i, X^\pi(S^i), W^{i+1}) = \left\{\begin{matrix}
D_{KL}(P^{i+1} \parallel U_{\mathcal{X}})-D_{KL}(P^{i} \parallel U_{\mathcal{X}}) & & \text{if} \ n \leq i \leq N-1, \\ 
0 & & \text{if} \ i=N .
\end{matrix}\right.
\end{align}
which is the increment in relative entropy between iteration $i+1$ and $i$.
Note that the data stream of this sequential decision process is $(S^0, x^0, W^1, S^1, x^1, W^2, ..., S^N, x^N)$. So with the transition function defined in section 3.2, we are able to calculate $P^{n+1}$ given $(S^n, X^\pi(S^n), W^{n+1})$. Then the surrogate value of a state $S^n$ at iteration $n \leq N$ is:
 \begin{align} \label{eq:30}
    \tilde{V}_n(S^n) 
    &= \mathbb{E}[\sum_{i=n}^{N} \tilde{C_i}(S^i, X^\pi(S^i), W^{i+1}) | S^n, x^n=(x_h,z)] \nonumber \\
    &= \mathbb{E}[D_{KL}(P^{N} \parallel U_{\mathcal{X}})-D_{KL}(P^{n} \parallel U_{\mathcal{X}}) | S^n, x^n=(x_h,z)]. 
\end{align}
At each iteration $n$, we choose to optimize the one-step lookahead value. Equivalently, set $\tilde{V}_{n+1}(S^{n+1}) = 0$ and obtain the one-step lookahead SBES policy:
\begin{align} \label{eq:31}
    X^{SBES}(S^n) &= \underset{h \in H^n, z \in \mathcal{X}}{\mathrm{argmax}} \Tilde{V}_n(S^n) \nonumber \\
                 &= \underset{h \in H^n, z \in \mathcal{X}}{\mathrm{argmax}} \mathbb{E}[\tilde{C}_n(S^n, x^n, W^{n+1}) | S^n, x^n=(h,z)] \nonumber \\
                 &= \underset{h \in H^n, z \in \mathcal{X}}{\mathrm{argmax}} \mathbb{E}[ D_{KL}(P^{n+1} \parallel U_{\mathcal{X}})-D_{KL}(P^{n} \parallel U_{\mathcal{X}}) ) | S^n, x^n=(h,z)] \nonumber \\
                 & = \underset{h \in H^n, z \in \mathcal{X}}{\mathrm{argmin}} \mathbb{E} [H(P^{n+1})-H(P^n) | S^n, x^n=(h,z)]
                 \coloneqq \nu^n (h,z).
\end{align}
The objective function $\nu^n$ in equation \eqref{eq:31} of the SBES algorithm is the same as the ES acquisition function despite the difference in sign. Another important merit of SBES is that with the probabilistic models stated in section \ref{model}, this objective function now has an analytical expression. Letting $(x_l,x_r)$ be the standard notation for indicating relative location of $(h, z)$, $\nu^n$ can be expressed as follows:
{\small\begin{align} \label{eq:sbes}
    \nu^n(h,z) &= \mathbb{E}[H(P^{n+1}) | S^n, x^n=(h,z)]-H(P^n) \nonumber \\
    &= [g^n(x_l,x_r)log_2(g^n(x_l,x_r)) + (1-g^n(x_l,x_r))log_2(1-g^n(x_l,x_r))](F^n(x_r)-F^n(x_l)-1) \nonumber\\
    &\quad - [\bar{g}^n(x_l,x_r)log_2(\bar{g}^n(x_l,x_r)) + (1-\bar{g}^n(x_l,x_r))log_2(1-\bar{g}^n(x_l,x_r))](F^n(x_r)-F^n(x_l)) \nonumber \\
    &\quad + U_1^n(x_l,x_r)\log_2(U_1^n(x_l,x_r)) + U_0^n(x_l,x_r)\log_2(U_0^n(x_l,x_r)).
\end{align}} Equation \eqref{eq:sbes} is the formula for a single-step lookahead. Yet, equation \eqref{eq:30} indicates that we are not limited to looking ahead just one-step; we can also perform multi-step lookaheads by maximizing the value function in \eqref{eq:30}.

\subsubsection{Optimizing the Surrogate Objective.} \label{discretize} 
We have derived the formula for the one-step lookahead entropy reduction objective above in equation \eqref{eq:sbes}. Even though this objective is more straightforward and easier to evaluate than the truth function $f$, it is still nonconvex in most scenarios. We now discuss how to optimize $\nu^n(h,z)$ when $\X$ is discrete and continuous respectively. Consider the following two possible structures of $\X$:

\begin{enumerate}
    \item \textbf{$\mathcal{X}$ has finitely many elements.} \\
    Note that the SBES objective function $\nu^n: H^n \bigtimes \mathcal{X} \rightarrow \mathbb{R}$ is a two-dimensional real-valued function. Both $H^n$ and $\mathcal{X}$ are finite in this case, so the complexity of finding the optimum is at most $\mathcal{O}(|\mathcal{X}|(N+1))$ (remember that $N$ is small). In other words, it is not hard to optimize $\nu^n$ over the domain even though it is nonconvex.
    \item \textbf{$\mathcal{X}$ is a compact interval in $\mathbb{R}$.}\\
    When $\mathcal{X}$ is a continuous and bounded set, optimizing $\nu^n(h,z)$ is equivalent to solving a deterministic one-dimensional nonconvex optimization problem over the feasible region $\mathcal{X}$ for $|H^n|$ times. So one solution is to use some well-studied deterministic, derivative-free global optimization methods, such as \textit{DIRECT} suggested in \citep{jones1993direct}, to optimize $\nu^n(h,z)$. 

    In spite of the availability of global optimization methods, we choose to discretize the search region $\mathcal{X}$ in order to approximate $\nu^n(h,z)$. To better facilitate the approximation, one recommended way of discretizing $\mathcal{X}$ is to draw $m$ samples from the posterior distribution $P^n$, as adopted by \citet{hernandez2014predictive}, \citet{shah2015parallel}, and \citet{wu2016parallel}. In general, the benefit of implementing this discretization method is it encourages exploration at the early stages of experimentation since we have less knowledge about $x^*$ and thus $P^n$ is more uniform over the entire search region. It also features exploitation later in the search as a more concentrated density $P^n$ will produce a collection of samples close to the estimated optimum. As a result, there are more chances at later iterations to check the local information around the estimated optimum.
    
    With the sampled belief model, another advantage of discretizing $\X$ based on $P^n$ is that $P^n$ is analytical and ready-for-use at iteration $n$. This avoids the relatively complicated sampling procedure involving a linear approximation of the GP prior using a feature map that \citet{hernandez2014predictive} introduces for optimizing the acquisition function \eqref{PES} of PES. 
    
    In the design of numerical experiments, the SBES algorithm generates $m$ points, denoted as $\mathcal{A}^n$, according to the up-to-date posterior distribution $P^n$ at every iteration and solve the following version of \eqref{eq:31}:
    \begin{align} \label{objective}
    X^{SBES}(S^n) &= \underset{h \in H^n, z \in \mathcal{A}^n}{\mathrm{argmin}} \nu^n(h,z).
    \end{align}
\end{enumerate}

Based on the discussion above, we argue that the SBES algorithm provides substantial simplification to implementation compared against other entropy-search based algorithms, in attribution to the fact that there is no extra step needed to evaluate the surrogate objective function \ref{objective}. In the case of discretization, it is also easier and more computationally tractable for SBES to sample the posterior distribution $P^n$ because the posterior distribution $P^n$ has numerical expressions. Consequentially, SBES is expected to improve on the computational time spent on optimizing the surrogate objective function.

\section{Theory} \label{theory} 
In this section, we provide an error bound for the final estimate $\bar{x}^N$, produced by the SBES algorithm, given that the truth function is in the set of sampled belief functions.

When noise is present (i.e. $\sigma > 0$), under the assumption that the sampled belief curves contains the underlying truth $f$, it is shown that the SBES algorithm is able to reduce the expected entropy of the posterior $P^n$, i.e. $\mathbb{E}[H(P^{n+1}) |x^n=(x,y),P^n]- H(P^n) \leq 0$. A single-period lower bound on the expected entropy reduction that is dependent on the function we are trying to learn is presented as well. 

For the following discussion, we use $x,y \in \mathcal{X}$ to refer to any two points in the search region $\mathcal{X}$.  When we want to emphasize the comparative location of $x,y$, we use $x_l$ and $x_r$ where $x_l = \mathrm{min} \{x,y\}$, $x_r = \mathrm{max} \{x,y\}$.

Before the formal proofs, we want to introduce two important assumptions, and then introduce the notations we use throughout the section.

\begin{assumption}\label{assume:true theta}
Suppose the underlying function $f$ is parameterized by $\theta$ and the true parameter is $\theta^*$, i.e. $f(x) = f(x |\theta^*)$. Then we assume that the set of belief parameters contains the true one : $\theta^* \in \Theta \coloneqq \{\theta_k\}_{k=1}^K$ ($|\Theta| = K$). Throughout the theory section, we suppose $\theta$ follows a distribution $p_{\theta}^n $ at iteration $n$ where $\theta_k \in \Theta$.
\end{assumption}

\begin{assumption}
Let $p_{\theta}^n = \{p^n_k\}_{k=1}^K$, where $p^n_k$ is defined in equation \eqref{def:prob_theta}, denote the posterior distribution of $\theta$ at iteration $n$ after observing the data $\mathcal{D}^n = (H^n, \hat{f}^n)$. Also denote the true distribution of $\theta$ to be $p_{\theta}^* \coloneqq \{\mathds{1}_{\{\theta_k=\theta^*\}}\}_{k=1}^K$.We assume that the prior distribution on $\theta$ satisfies $p^0_{k} > 0$.
\end{assumption}

Recall that $\hat{y}^{n+1} \coloneqq \mathds{1}_{\{\hat{f}(x_l^n) \leq \hat{f}(x_r^n)\}} $. At any state $S^n$, the tuple $(x_l^n, x_r^n)$ is determined by a policy $\pi$ that picks two points $X^{\pi}(S^n) = (x_h^n, z^n)$. The binary variable $\hat{y}^{n+1}$ is then a function of $(\hat{f}(x_h^n), f_{\theta}(z^n), W^{n+1})$. Note that $f_{\theta}(z^n)$ is the noiseless function value evaluated at $z^n$ given that the true parameter is $\theta$. In other words, given a decision $X^{\pi}(S^n)$ the distribution of $\hat{y}^{n+1}$ is determined jointly by the distribution of $\theta$ and $W^{n+1}$.
In the following sections, we use $\hat{y}^{n+1}(X^\pi(S^n))$ to indicate the binary random variable $\hat{y}^{n+1}$ under a policy $\pi$ at the state $S^n$ whenever we want to emphasize the relation between $\hat{y}^{n+1}$ and $\pi$. When the policy $\pi$ is not clearly defined, we use $\hat{y}^{n+1}(x^n)$, in which $x^n$ is some decision made at iteration $n$.

To produce a good estimate of the location of the maximizer $x^*$, the learning algorithm should collect as much information as possible throughout $N$ experiments. The information gain about $x^*$ is quantified in the mutual information between the random variable $X^*$ and the $N$ random observations: $I(X^*; \{\hat{y}_{\pi}^{n}\}_{n=1}^N|S^n)$. In the following, we investigate the difference between the information gain of the SBES decision and the information gain of the conceptually optimal policy $\pi^*$, given the underlying true parameter $\theta^*$. This section starts with specifying the difference between predictive mutual information and the perfect mutual information, as well as the general assumptions we make. Then, we proceed to prove an error bound on the one-step perfect mutual information generated by the SBES policy. It is shown that this error bound is tied to the quality of our estimation of the true parameter $\theta^*$.

\begin{assumption}[Distinct Belief Optimizers]
All of the sampled belief functions have distinct optimizers. Define $x^*_k \coloneqq \text{argmax} \ f(x | \theta_k)$. Then $x^*_i \neq x^*_j, \ \forall i \neq j \in \{1, ..., K \}$. 
\end{assumption}

\begin{assumption}[Existence of One-to-one Mapping] \label{assume_4}
Let a general probability space of $\theta$ be $(\Theta, \mathcal{F}, \Prob_{\theta})$. It is assumed that there exists a quantizer $q:\mathcal{X} \rightarrow \{1, ..., K \}$ and a partition $\mathcal{P}_q \coloneqq \{\mathcal{X}_1,...,\mathcal{X}_K \}$ associated with $q$ such that $x^*_k \in \mathcal{X}_k, \ \forall k=1, ..., K$. Let $\sigma(\mathcal{P}_q)$ be the sigma algebra generated by $\mathcal{P}_q$. Then the random variable $X^*: \Theta \rightarrow \mathcal{P}_q$ is well-defined with the induced probability space $(\mathcal{P}_q, \sigma(\mathcal{P}_q), \Prob_{X^*})$.
\end{assumption}

\begin{remark}
Assumption 4 is saying that there exists a one-to-one mapping between the event $\{\theta = \theta_k\}$ and $\{X^* \in \mathcal{X}_k\}, \  \forall k = 1,..., K$. In other words, knowing any of the events directly implies the other event.
\end{remark}

\begin{definition}[Predictive Mutual Information]
Let $p(\hat{y}^{n+1} | \theta)$ denote the distribution of $\hat{y}^{n+1}$ parameterized by $\theta$. The predictive mutual information $\hat{I}$ at a state $S^n$ is defined as:
\begin{align}
    \hat{I}(X^*;\hat{y}^{n+1}|S^n) 
    & \coloneqq H(p(\hat{y}^{n+1} | S^n)) - \E_{X^*|S^n}[H(p(\hat{y}^{n+1} | S^n, X^*))] \nonumber \\
    & \coloneqq  H(p(X^*| S^n)) -\mathbb{E}_{\hat{y}^{n+1} |S^n \sim p(\hat{y}^{n+1} | \theta)}[H(p(X^*|S^n, \hat{y}^{n+1}))] \nonumber \\
    & = H(p(X^* | S^n)) -\mathbb{E}_{\theta}[\mathbb{E}_{\hat{y}^{n+1} |S^n, \theta}[H(p(X^*|S^n, \hat{y}^{n+1}))]]. 
\end{align}
Note that $S^n$ encodes the posterior distribution of $\theta$ given all the information up to iteration $n$, i.e. $p^n_{\theta}$. This implies that as we update the posterior distribution of $\theta$, the distribution of $\hat{y}$: $p(\hat{y}^{n+1} | \theta)$ also changes.
\end{definition}

\begin{definition}[Perfect Mutual Information]
If we were given perfect information at the state $S^n$, i.e. $\theta = \theta^*$ with probability 1, then we know the true distribution of $\hat{y}^{n+1}$: $p(\hat{y}^{n+1} | \theta^*)$ as well as $p(X^* | \theta^*)$. Define the perfect mutual information as:
\begin{align}
    I^*(X^*;\hat{y}^{n+1}|S^n) & = 
    H(p(X^* | S^n)) - \mathbb{E}_{\hat{y}^{n+1} |S^n \sim p(\hat{y}^{n+1} | \theta^*)}[H(p(X^*|S^n, \hat{y}^{n+1}))]]   \nonumber \\
    & = H(p(\hat{y}^{n+1} | S^n)) - \E_{X^*|S^n,\theta^*}[H(p(\hat{y}^{n+1} | S^n, X^*))].
\end{align}
\end{definition}

\begin{lemma} \label{max_pmi}
For $\hat{y}^{n+1}(X^\pi(S^n))$ produced by any policy $\pi$, its predictive mutual information is no greater than the predictive mutual information of $\hat{y}^{n+1}(X^{SBES}(S^n))$:
\begin{align}
    \hat{I}(X^*;\hat{y}^{n+1}(X^\pi(S^n))|S^n) \leq \hat{I}(X^*;\hat{y}^{n+1}(X^{SBES}(S^n))|S^n).
\end{align}
\end{lemma}
\proof{Proof.}
Recall that the decision of SBES $X^{SBES}(S^n)$ is determined by maximizing the predictive entropy reduction of the posterior distribution of $X^*$ with the pair $(x_h, z) \in H^n \times \mathcal{X}$ under the parameter $\theta$:
\begin{align}
    X^{SBES}(S^n) 
    = & \underset{x_h \in H^n, z \in \mathcal{X}}{\mathrm{argmin}} \mathbb{E}_{\theta, W^{n+1}} [H(P^{n+1})-H(P^n) | S^n, x^n=(x_h,z)] \nonumber \\
    =& \ \underset{x_h \in H^n, z \in \mathcal{X}}{\mathrm{argmax}} \ H(p(X^* | S^n)) - \mathbb{E}_{\hat{y}^{n+1}}[H(p(X^* |S^n, \hat{y}^{n+1}(x^n))] \nonumber \\
    =& \ \underset{x_h \in H^n, z \in \mathcal{X}}{\mathrm{argmax}} \ H(p(X^* | S^n)) - \mathbb{E}_{\theta}[\mathbb{E}_{\hat{y}^{n+1}|\theta}[H(p(X^* |S^n, \hat{y}^{n+1}(x^n)]] \nonumber \\
    =& \ \underset{\pi}{\mathrm{argmax}} \ \hat{I}(X^*;\hat{y}^{n+1}(X^\pi(S^n))|S^n). 
\end{align}
\endproof

In the next theorem, we bound the error between the predictive mutual information produced by the SBES algorithm and the maximum perfect mutual information. The idea behind the perfect mutual information is that the quality of the next query point we choose is measured under the true parameter $\theta^*$. Ideally, we should pick a point that maximizes the mutual information between the true optimum $x^*$ and the next observation. However, the reality is that we do not know this true parameter $\theta^*$, so we use the predictive mutual information as an estimate of the perfect mutual information, which produces an error. The following theorem provide a bound on this error under the measure of perfect mutual information.

\begin{theorem} \label{bound_error}
Let $\pi^*$ be the optimal policy that maximizes the perfect mutual information. The error of a single-step information gain between $\pi^*$ and the SBES algorithm is bounded by the KL-divergence from the posterior distribution of $\theta$ to the point mass measure $p^*_{\theta} \coloneqq \{\mathds{1}_{\{\theta_k=\theta^*\}}\}_{k=1}^K$. That is:
\begin{equation}
    \begin{aligned}
    & \left| I^*(X^*;\hat{y}^{n+1}(X^{\pi^*}(S^n))|S^n) - I^*(X^*;\hat{y}^{n+1}(X^{SBES}(S^n))|S^n) \right| \nonumber \\
     \leq & ~ 4\mathbb{P}(\theta^n \neq \theta^*) \nonumber \\
    \leq & ~ \mathcal{O} \left( \sqrt{D_{KL}(p^*_{\theta} || p^n_{\theta})} \right).
    \end{aligned}
\end{equation}
For the sake of simplicity, we use $\hat{y}_{\pi^*}$ and $\hat{y}_{SBES}$ to denote $\hat{y}(X^{\pi^*}(S^n))$ and $\hat{y}(X^{SBES}(S^n))$ respectively.
\end{theorem}

\proof{Proof.}
\begin{align}
    & I^*(X^*;\hat{y}_{\pi^*}^{n+1}|S^n) - I^*(X^*;\hat{y}_{SBES}^{n+1}|S^n) \nonumber \\
   =& I^*(X^*;\hat{y}_{\pi^*}^{n+1}|S^n) - \hat{I}(X^*;\hat{y}_{\pi^*}^{n+1}|S^n) +  \hat{I}(X^*;\hat{y}_{\pi^*}^{n+1}|S^n) - \hat{I}(X^*;\hat{y}_{SBES}^{n+1}|S^n) \nonumber \\
   & + \hat{I}(X^*;\hat{y}_{SBES}^{n+1}|S^n) - I^*(X^*;\hat{y}_{SBES}^{n+1}|S^n) \nonumber \\
   \leq & I^*(X^*;\hat{y}_{\pi^*}^{n+1}|S^n) - \hat{I}(X^*;\hat{y}_{\pi^*}^{n+1}|S^n) + \hat{I}(X^*;\hat{y}_{SBES}^{n+1}|S^n) - I^*(X^*;\hat{y}_{SBES}^{n+1}|S^n).
\end{align}
where the last inequality follows from lemma \ref{max_pmi}.
\\
For any policy $\pi$, we have:
\begin{align}
    & \left|I^*(X^*;\hat{y}_{\pi}^{n+1}|S^n) - \hat{I}(X^*;\hat{y}_{\pi}^{n+1}|S^n)\right| \nonumber \\
    = & \left|\E_{X^*|S^n}[H(p(\hat{y}_{\pi}^{n+1} | S^n, X^*))] -\E_{X^*|S^n, \theta^*}[H(p(\hat{y}_{\pi}^{n+1} | S^n, X^*))] \right| \nonumber \\
    = & \left|\int_{\mathcal{X}} p(X^* \in dx|S^n) H(p(\hat{y}_{\pi}^{n+1} | S^n, X^* \in dx)) - \int_{\mathcal{X}} p(X^* \in dx|S^n, \theta^*) H(p(\hat{y}_{\pi}^{n+1} | S^n, X^* \in dx)) \right| \label{measure_x}\\ 
    = & \left|\int_{\Theta} p(\theta^n \in d\theta|S^n) H(p(\hat{y}_{\pi}^{n+1} | S^n, X^*(d\theta))) - H(p(\hat{y}_{\pi}^{n+1} | S^n, X^*(\theta^*)) \right| \label{measure_theta}\\
    = & \left|\sum_{k=1}^K p^n_k H(p(\hat{y}_{\pi}^{n+1} | S^n, \theta_k)) -  H(p(\hat{y}_{\pi}^{n+1} | S^n, \theta^*) \right| \label{replace_x}  \\
    = & \left|\sum_{k=1}^K (p^n_k-\mathds{1}_{\{\theta_k=\theta^*\}}) H(p(\hat{y}_{\pi}^{n+1} | S^n, \theta_k)) \right|  \nonumber \\
    \leq & \norm{p^n_\theta - p^*_\theta}_1 \norm{H(p(\hat{y}_{\pi}^{n+1} | S^n, \theta))}_\infty   \nonumber \\
    \leq & \norm{p^n_\theta - p^*_\theta}_1 = 2\mathbb{P}(\theta \neq \theta^*).
    \label{holder}
\end{align}
The equality between \eqref{measure_x} and \eqref{measure_theta} follows from assumption \ref{assume_4}, whereas the replacement of $X^*$ with $\theta$ in $\eqref{replace_x}$ occurs because when $X^* \in \mathcal{X}_k$, where $\mathcal{X}_k$ is the partition subset which corresponds to $\theta^k$, the conditional distribution of $\hat{y}^{n+1}$ is the same as the one conditioned on $\theta^n=\theta^k$. In other words, $X^*$ influences the distribution of $\hat{y}^{n+1}$ through its association with the underlying function, which is a piece of information that $\theta^n$ also encodes. Inequality \eqref{holder} is immediate from the fact that $\hat{y}^{n+1}$ is a binary variable whose entropy has maximum equal to $1$.

The proof is finished by Pinsker's inequality:
\begin{align}
 & \left| I^*(X^*;\hat{y}_{\pi^*}^{n+1}|S^n) - I^*(X^*;\hat{y}_{SBES}^{n+1} |S^n) \right| \nonumber\\
 \leq & \left|I^*(X^*;\hat{y}_{\pi^*}^{n+1}|S^n) - \hat{I}(X^*;\hat{y}_{\pi^*}^{n+1}|S^n) \right| + \left| \hat{I}(X^*;\hat{y}_{SBES}^{n+1}|S^n) - I^*(X^*;\hat{y}_{SBES}^{n+1}|S^n) \right| \nonumber\\
    \leq & 2\norm{p^n_\theta - p^*_\theta}_1 = 4\norm{p^n_\theta - p^*_\theta}_{\textit{TV}} \nonumber\\
    \leq & 4\sqrt{2D_{KL}(p^*_{\theta} || p^n_{\theta})}. \nonumber
\end{align}
\endproof

The take-away from theorem \ref{bound_error} is that at every state, the error between the information gain of the SBES algorithm and the information gain of the optimal policy is bounded by the estimation error of $\theta^n$, which is a model dependent quantity.

\begin{corollary}(Lower Bound on One-step Mutual Information)
\begin{align}
    \mathrm{max} \left\{ I^*(X^*;\hat{y}_{\pi^*}^{n+1}|S^n) - 4\sqrt{2D_{KL}(p^*_{\theta} || p^n_{\theta})}, \ 0\ \right\} \leq I^*(X^*;\hat{y}_{SBES}^{n+1}|S^n).
\end{align}
\end{corollary}
\proof
Given that $I(X;Y) > 0$ if $X,Y$ are not independent \citep{John2019}, the statement holds directly from theorem \ref{bound_error}.
\endproof

\section{Numerical Experiments} 
\label{section5:experiment}
We present the results of two types of experiments. First, the SBES algorithm is tested on maximizing synthetic functions, both parameterized and black-box, in comparison to other well-studied algorithms. In all experiments, we implement the SBES algorithm under different noise levels as well as multiple initializations and illustrate that it outperforms other algorithms especially at a medium noise level. To further demonstrate the power of SBES, it is employed as the stepsize rule of the Gradient Descent (GD) algorithm, which is applied to finding the optimum of various multidimensional functions. It is observed that SBES is more robust than prevalent stepsize rules such as RMSProp and AdaGrad to the choice of initialization. Moreover, GD converges much faster with SBES as a parametric stepsize rule.

\subsection{Synthetic Functions} 
Two types of synthetic functions are investigated in testing the performance of the SBES algorithm: unimodal parametric functions (Gaussian, Gamma and Beta) and well-known benchmark functions (Mccormick and Ackley) noted in \cite{synthetic}. For the first category, the SBES algorithm is provided with the exact parametric family and a set of parameters that include the true one. In addition, there are two versions of the SBES algorithm. One of them is provided with the real scale of the test function, while the other version, called SCALE-SBES, has to learn the true scale as the experiment proceeds. For the second category of test functions, since we don't know which family the test function belongs to, we pick several parametric families as an approximation of the test function and run the SCALE-SBES algorithm based on the parametric beliefs of our choice. The selected parametric families are specified in \cref{tab:1}.

In terms of the benchmark algorithms against which the SBES algorithm competes, we choose the response surface method (RSM), GP-UCB, GP-EI, PES and MES where the GP kernel’s parameters are optimized along the way. Each experimental result is computed with $15$ randomly chosen initializations from a Latin hypercube design and $20$ realizations of each initialization. Different noise levels are also considered when comparing the performances of different algorithms. Since the scale of test functions can be different, to make noise levels more comparable we use noise ratio $(\gamma)$, defined as the ratio between standard deviation of the noise in observations and the maximum difference in underlying function values:
\begin{equation} \label{eq:noise_ratio}
    \gamma \coloneqq \frac{\sigma}{|f_{\text{max}} -f_{\text{min}}|}.
\end{equation}
For example, the 2-dimensional Rosenbrock function has $|f_{\text{max}} -f_{\text{min}}| = 1000$, 
so a noise ratio equal to $0.5$ means $\sigma = 500$.

\Cref{tab:1} shows the immediate regret of all algorithms in log scale under different noise levels, in which a smaller regret means a more accurate estimation. The results illustrate that SBES is the absolute winner of the in-model experiments with the outperforming gap maximized at a medium noise level. Despite the effort to learning the true scale, SCALE-SBES has a similar performance to SBES, suggesting that SBES's performance is robust against the variation in scale as long as the true parameter is in its sampled belief set. This is a significant assumption as in reality, people usually have limited information about the magnitude of the underlying function aside from which family of parametric functions it resembles.
\begin{table}[h!]
\centering
\small
\caption{\label{tab:1}  Performance on standard test functions.}
\renewcommand{\arraystretch}{1.2}
\subcaption{Results on comparing benchmark algorithms , SBES and SCALE-SBES on parametric functions.}
\smallskip
{\begin{tabular}
{c m{1cm}|c c c c c c c}
\hline
\text{Test Function}                           & \text{Noise} \text{Level}  & GP-EI & GP-UCB & RSM & PES & MES & SBES & SCALE-SBES \\\hline
\multirow{3}{*}{Gamma}      &\text{low}  &  -2.88  &  -2.79   &  1.17   &  -0.06   &   -3.63  & \textbf{-6.33}     & -6.28\\
                            & \text{mid}  & -0.31   &  -0.25   &   1.16  &   0.02  &  -1.72   &   \textbf{-3.98}    &    -3.52        \\
\text{$k=9$,$\lambda=1$}    &\text{high}&  -0.07  & 0.01    &   -0.25  &  0.15   &  \textbf{-0.58}   &   -0.54   &  -0.43          \\\hline
\multirow{3}{*}{Beta}       &\text{low}   &   0.25 &    0.18 &  1.07   &  0.15   &  -2.88   &  \textbf{-4.31}    &   -4.22         \\ 
                            & \text{mid} &  0.36  &  0.64   &    0.91 &   0.18  &   -1.27  &   \textbf{-4.26}   &  -3.58          \\
\text{$\alpha=3$,$\beta=18$}   & \text{high} & 0.64   &    0.32 & 1.08    & 0.14    &  -0.02   &   \textbf{-1.06}   &     -0.34       \\\hline
\multirow{3}{*}{Gaussian}   &\text{low}   &  -3.57  &  -3.65   &    1.02 &   0.08  &  -3.78   &   \textbf{-5.88}   &  -5.84         \\ 
                            & \text{mid} &  -2.02  & -2.21    &  1.10   &  0.10   &  -2.08   &    \textbf{-5.61}  &    -5.06        \\
\text{$\mu=7.5$,$\sigma=1$} &\text{high}  & -0.08   &    -0.30 & 0.92    &  0.13   &  -0.68   &   \textbf{-1.06}   &    -0.75        \\\hline
\multicolumn{2}{c|}{\text{overall average}}   & -0.85   &    -0.89 & 0.91    &  0.099   &  -1.85   &   \textbf{-3.67}   &    -3.34        \\\hline
\end{tabular}
\bigskip
\subcaption{Results on comparing benchmark algorithms and SCALE-SBES on Mccormick: $f(x) = -sin(x)-x^2+1.5x+10$, Ackley: $f(x)= 4e^{-\|x\|}+e^{cos(x)}-4-e$.}
\smallskip
\begin{tabular}
{c m{1cm}|c c c c c c c}
\hline
\text{Test Function} & Noise Level  & GP-EI & GP-UCB & RSM & PES & MES                          & \multicolumn{2}{p{3cm}}{Parametric Family Used by SCALE-SBES}  \\
\hline
\multicolumn{7}{l}{} & \multicolumn{1}{l}{Gamma} & \multicolumn{1}{l}{Quadratic} \\\hline

\multirow{3}{*}{Mccormick} &\text{low}   &  -2.90  &  -3.49   &  -2.20   &  -0.08   &  \textbf{-3.36}   & -0.69    &  -2.50          \\
           & \text{mid} & -1.49   & -1.82    & -1.06 &   -0.19  & -1.88    &    -0.80  & \textbf{-2.17}           \\
                            & \text{high} &  -0.25  & \textbf{ -0.99} &  -0.69   & -0.08    &  -0.97   &    -0.43  & -0.62         \\\hline
\multicolumn{1}{l}{}     & \multicolumn{6}{l}{}                                                                                                                  & \multicolumn{1}{l}{Gaussian} & \multicolumn{1}{l}{Quadratic} \\
\hline
\multirow{3}{*}{Ackley}     & \text{low}  &  0.74  &    0.74 &   -0.27  &  -0.13   &   -0.27  &   \textbf{-2.21}   & -0.63          \\
                            & \text{mid} &  0.74  &  0.74   &   -0.5  &  -0.09   &   -1.01  &   \textbf{-1.42}   &   -1.32   \\
                            & \text{high} & 0.74   &    0.72 &   -0.22  &  -0.05   &  \textbf{-0.34}   &   0.24   &  -0.26        
                           \\ \hline
\multicolumn{2}{c|}{\text{overall average}}   & -0.40   &    -0.68 & -0.82    &  -0.10   &  \textbf{-1.31}   &   -0.89   &    -1.25       \\\hline
\end{tabular}}
\\
\medskip
{\textit{Notes.} Each row summarizes the mean immediate regret in log10 scale from 900 runs after 30 iterations of each policy at the following corresponding noise level. Low noise: $\gamma \in [0.003, 0.007]$; mid noise: $\gamma \in [0.03, 0.125]$;  high noise: $\gamma \in [0.3, 0.5]$, where $\gamma$ is defined in \cref{eq:noise_ratio}. The lowest regret of each row is highlighted in bold.}
\end{table}

When testing on black-box benchmark functions (Mccormick and Ackley), SBES still beats all of the competing algorithms in the presence of medium noise as shown in table \ref{tab:1}. An interesting observation is that SBES has a better performance on Ackley than on Mccormick, while all the GP-based methods are more successsful on Mccormick. The difference is that Mccormick features a smooth and convex curvature whereas Ackley is a nonconvex, triangle-shaped function with a sharp and non-differentiable optimum. This phenomenon implies that SBES is competent even when the assumptions about the underlying function are relaxed. It is worth mentioning that MES is the best among all the competing algorithms in terms of robustness, especially when the underlying function is smooth. Yet SBES either outperforms or is competitive against MES across all test functions.

\subsection{Stepsize Experiments} 
An important application of SBES is to conduct the one-dimensional line search for (stochastic) Gradient Descent, which is an essential algorithm to modern machine learning methods. Finding the stepsize is typically handled using a parametric stepsize rule such as Adagrad or RMSProp which has to be tuned to capture not only the characteristics of the problem, but also the choice of starting point. However, practitioners often find that popular stepsize formulas suffer from bad choices of the starting point of the GD algorithm. An algorithm that is less sensitive to the choice of starting point would represent a major contribution to stochastic gradient descent. We use this section to demonstrate that SBES can exactly achieve the described objective.

\begin{figure}[t]
    \centering
    \caption{Average reduction in the initial distance to the global maximum within 10 iterations on far starting points.}
    \includegraphics[width=\linewidth]{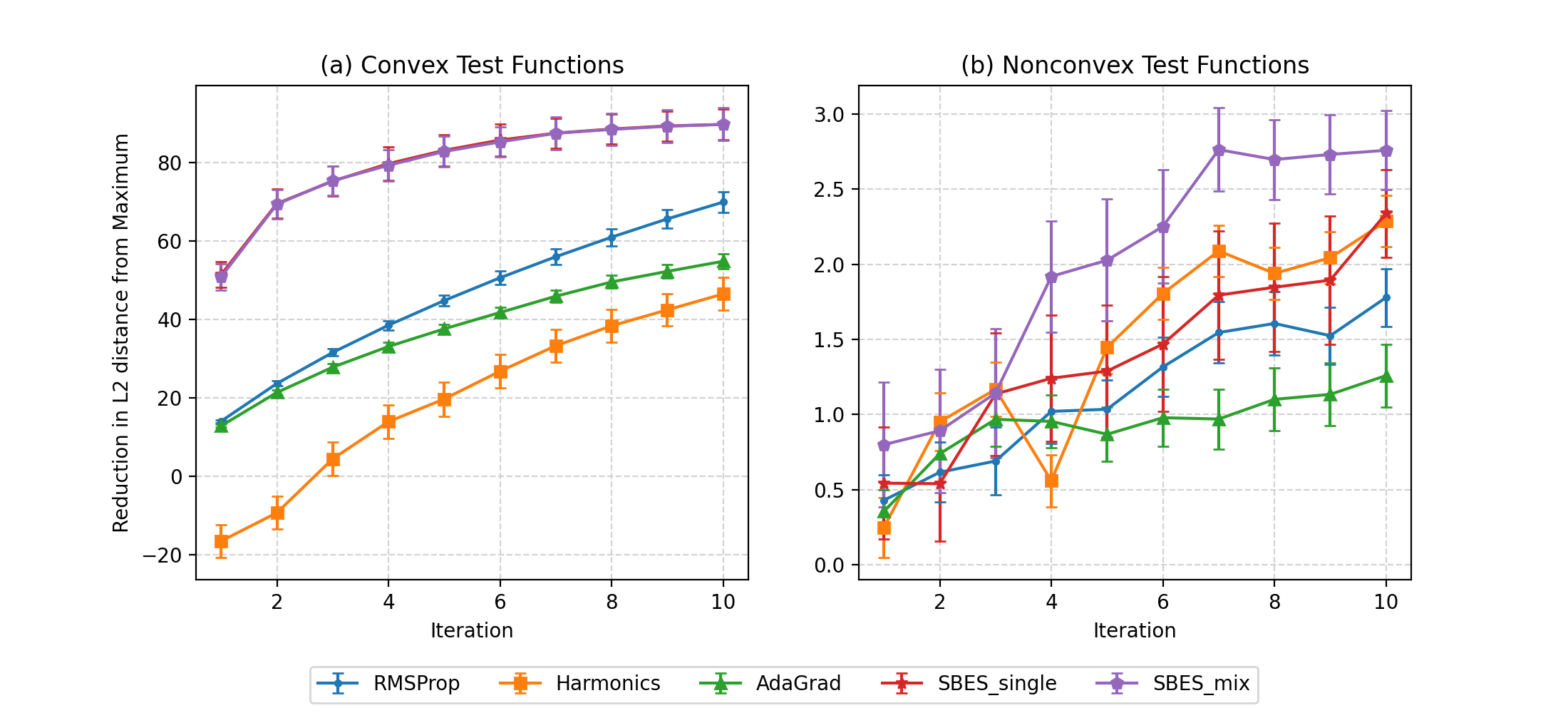}
    \label{fig:far_compare}
    \small
    {\textit{Note.} For all stepsize experiments, SBES runs 5 iterations to produce one single stepsize.}
\end{figure}

In the following experiments, we use stochastic Gradient Descent(SGD) to find the global optimum of multidimensional unimodal functions defined on compact hypercubes. The noise present in the evaluation of the truth functions has standard deviation of $0.1$. All stochastic gradients are estimated via finite-difference stochastic approximation (FDSA).The SBES algorithm is applied as a stepsize rule of SGD, which is compared against a selected subset of ad-hoc benchmark stepsize formulas including harmonic, AdaGrad and RMSProp. Our assumption is still that the underlying function is noisy and expensive to evaluate, resulting in hard-to-obtain gradients. For this reason, we restrict our attention to the performance of each stepsize rule after 10 iterations of SGD. We use the distance between the furthest vertex and the known optimum $x^*$ as a reference for picking the starting point of SGD. Denote this distance to be $d$. Then an initial point is randomly chosen from a region in-between two spheres defined by $\mathcal{B}(x^*,r_2) \backslash \mathcal{B}(x^*,r_1)$, in which $r_1, r_2$ are a proportion of $d$. See \cref{tab:2} for the specific choice of $r_1, r_2$ for each category of initial points.

There are two categories of unimodal test functions we have used: convex and non-convex. For the convex test functions, we picked the Bohachevsky, Rotated Hyper-Ellipsoid and Sum of Different Powers functions. The Bohachevsky function is defined to be 2-dimensional while the dimension of the other two functions can be varied. We run the simulations with the choices of dimension for the other two test functions to be 5, 10 and 20. The non-convex test functions consist of multivariate Gaussian density functions in the hypercube $x_i \in [0,5]$ for $i=2, 10$ with different means and covariance matrices.

\begin{table}[]
\small
\caption{Average reduction in L2-Distance from the global optimum ($\norm{x^0-x^*}$) after 10 iterations of stochastic Gradient Descent with different starting points.}
\label{tab:2}
\renewcommand{\arraystretch}{1.2}
\centering
\begin{tabular}
{c c| c c c c c}
\hline
& \text{initialization} & Harmonic & RMSProp & AdaGrad & SBES-single & SBES-mix \\ \hline
\multirow{3}{*}{\begin{tabular}[c]{@{}l@{}}nonconvex\\ unimodal\\ functions\end{tabular}} & close            &     -2.81    &    -2.06     &    -1.50     &       \textbf{0.36}        &     0.31            \\  
  & medium           &    0.38       &    0.30     &    0.26     &    1.48           &       \textbf{1.53}          \\ 
  & far              &    2.29  &  1.78   &    1.26     &    2.34     &  \textbf{2.76}                             \\ \hline
\multirow{3}{*}{\begin{tabular}[c]{@{}l@{}}convex\\ unimodal\\ functions\end{tabular}}    & close            &  -21.87         &    20.07     &     21.74    &       \textbf{22.54}        &      22.52           \\  & medium           &     23.48      & 56.14        &    45.28     &      66.97         &  \textbf{67.12}               \\  & far              &    46.49       & 69.89        &     54.81    &      89.67         &  \textbf{89.70}               \\ \hline
\multicolumn{2}{c|}{\text{overall average}}   & 7.99   &    24.35 & 20.31    &  30.56   &  \textbf{30.66}  \\ \hline
\end{tabular} \\
\medskip
{\textit{Notes.} Each test function is tested on 20 random starting points from each category. The radii of spheres chosen for each category of initial points are: close ($r_1=0, r_2=\frac{1}{4}d$), medium ($r_1=\frac{1}{4}d, r_2=\frac{3}{4}d$), far ($r_1=\frac{3}{4}d, r_2=d$). The largest reduction in $\norm{x^0-x^*}$ of each row is highlighted in bold.}
\renewcommand{\arraystretch}{1.5}

\end{table}

We provide two versions of the SBES stepsize rule: SBES-single and SBES-mix. The sampled beliefs in SBES-single all have the same shape and only vary in the horizontal localities. For example, we propose a family of gaussian density curves with the same variance but different means in the non-convex experiments. The second version SBES-mix employs sampled beliefs that have different localities as well as shapes, which can be a result of variation in parameters within the same parametric family or different choices of parametric family. We can see from \cref{tab:2} that both SBES-single and SBES-mix achieve the highest reduction in the distance between a Gradient Descent iterate and the global optimum in all categories of the starting position, especially in the sector of far initialization (highlighted in bold).  Figure \ref{fig:far_compare} further stresses that the L2-distance from the starting point to the optimum has been substantially reduced only after the first iteration.

\begin{figure}
    \centering
    \caption{Comparison of Steps for Gradient Descent on Maximizing a 2-dimensional Gaussian Density Function.}
    \includegraphics[width=\linewidth]{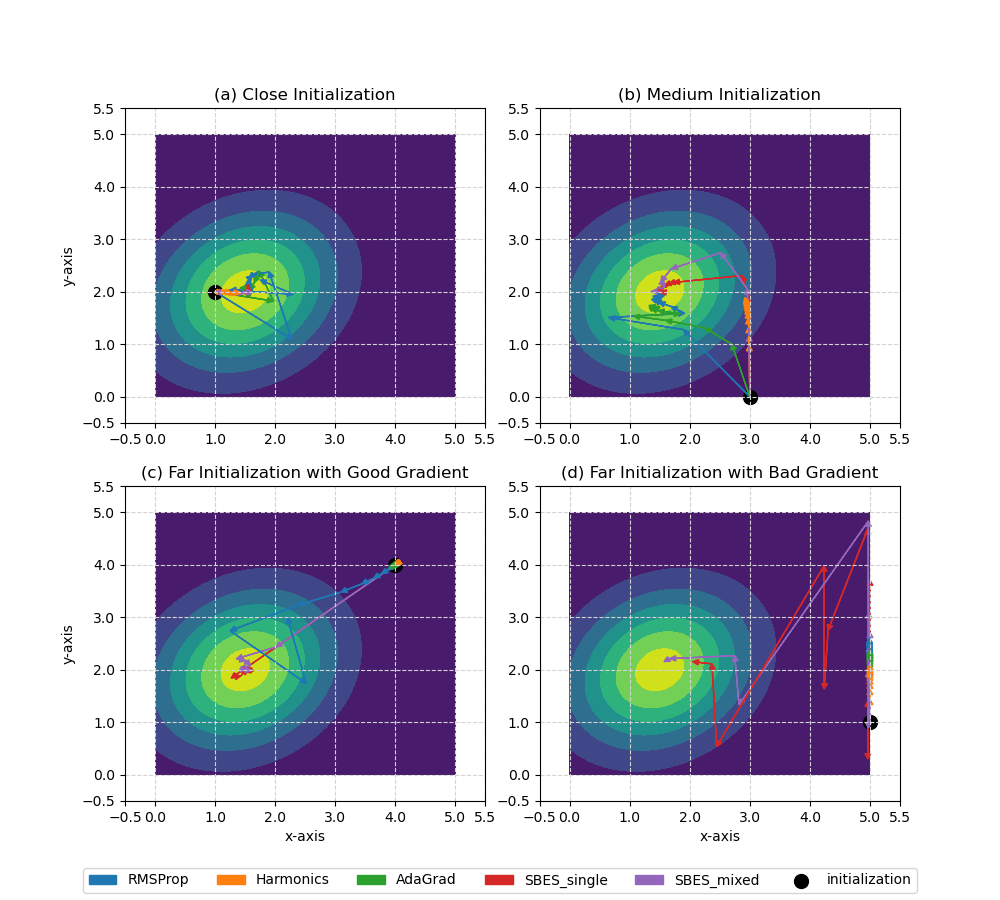}
    \small
    {\textit{Note.} To better demonstrate the behaviors of different stepsize algorithms, the results shown in the above heatmaps are based on noiseless gradients.}
    \label{fig:heatmap}
\end{figure}

One key note is that all of harmonic, RMSProp and AdaGrad need to be tuned for different underlying functions in order to perform reasonably. But tuning is not a remedy for the problem of initialization as well. It is common that the pre-tuned stepsize rules only converge on some initial points but fail on many others. In practice, retuning a stepsize formula for different starting points on the same problem is simply unrealistic. To the opposite, SBES uses the exact same set of parametric beliefs within a category (eg. convex) of functions and yields consistent performance across various types of initialization. If we refer to \cref{fig:heatmap}, it is easy to see that RMSProps and AdaGrad zig-zag around the optimum and harmonic is too conservative in its steps. Meanwhile, SBES chooses a stepsize based on how close the current iterate is to the optimum.

\section{Conclusion}
We have proposed a novel Entropy Search (ES) based algorithm, Sampled Belief Entropy Search (SBES), built upon a discrete parametric belief model instead of the Gaussian Process (GP) regression. While SBES inherits the logic of maximizing the one-step expected reduction in the entropy of the distribution of $x^*$, it successfully transforms the computationally intractable objective function of ES into an analytic formula under the assumption that the truth function is one-dimensional and unimodal. This is a significant step in the effort to simplify the computation of expected entropy reduction, as are attempted in Predictive Entropy Search (PES) and Max-Value Entropy Search (MES). While both PES and MES still require numerical approximation as well as sampling techniques, the surrogate objective function of SBES is ready-to-use. 

In addition to the computational benefits, SBES achieves the following accomplishments: 1) We have proved a lower bound on the one-step information gain produced by SBES, which occurs to be a problem-dependent quantity. 2) Experiments with synthetic functions also show that SBES often outperforms both ES based algorithms and non-ES based algorithms within a limited amount of experimental budgets. 3) Among all the test functions we use, an important observation is that SBES is able to handle functions with a non-smooth structure much better than algorithms modeled on GP, such as PES, MES, GP-UCB and GP-EI. The reason could be that the quality of the estimation produced by the GP based methods depends on how well the Gaussian Process regression captures the curvature of the truth function. Hence, in the scenario that the truth function is non-smooth or noisy, GP falls short of reproducing the underlying function precisely, thus leading to inaccurate predictions. In contrast, SBES focuses on learning the right comparative relations between any two points through the probability of correct region assignment. 4) On top of this, SBES makes estimates according to the posterior distribution of $x^*$, which in combination with the previous feature enables SBES to combat noisy and expensive experiments. 5) Lastly, in the experiments where SBES is applied as a stepsize rule for Gradient Descent (GD), it is evident that SBES is less sensitive to the choice of the starting point of GD than ad hoc stepsize formulas such as harmonic, RMSProp and AdaGrad. SBES is also more robust in the sense that it does not need much tuning across different types of functions to produce good results, unlike harmonic, RMSProp and AdaGrad. 

%
%
%


\bibliographystyle{apalike}
\bibliography{myrefs} 


\end{document}